\documentclass[12pt]{amsart}
\usepackage{amsfonts, amsmath, latexsym, epsfig}
\usepackage{amssymb}
\usepackage{epsf}
\usepackage{url}

\newcommand{\RR}{\ensuremath{\mathbb{R}}}

\newcommand{\QQ}{\ensuremath{\mathbb{Q}}}

\newcommand{\ZZ}{\ensuremath{\mathbb{Z}}}
\newcommand{\PP}{\ensuremath{\mathbb{P}}}

\newtheorem{theorem}{Theorem}
\newtheorem{corollary}{Corollary}
\newtheorem{lemma}{Lemma}

\def\QuotS#1#2{\leavevmode\kern-.0em\raise.2ex\hbox{$#1$}\kern-.1em/\kern-.1em\lower.25ex\hbox{$#2$}}

\DeclareMathOperator{\Sym}{Sym}

\DeclareMathOperator{\GL}{GL}

\newcommand{\op}{\operatorname}

\newcommand{\Sat}{{\mathcal A_g^{\op {Sat}}}}
\newcommand{\Vor}[1][g]{{\mathcal A_{#1}^{\op {Vor}}}}
\newcommand{\Perf}{{\mathcal A_g^{\op {Perf}}}}

\newcommand{\Igu}{{\mathcal A_g^{\op {Igu}}}}
\newcommand{\tor}{{\mathcal A_g^{\op {tor}}}}
\newcommand{\torb}{{\mathcal A_{g,{\op{smooth}}}^{\op {tor}}}}
\newcommand{\tors}{{\mathcal A_{g,{\op{simp}}}^{\op {tor}}}}
\newcommand{\torl}{{\mathcal A_g^{\op {tor}}(\ell)}}

\begin{document}

\author[M. D. Sikiri\'c]{Mathieu Dutour Sikiri\'c}
\address{Mathieu Dutour Sikiri\'c, Rudjer Boskovi\'c Institute, Bijenicka 54, 10000 Zagreb, Croatia}
\email{mathieu.dutour@gmail.com}

\author[K. Hulek]{Klaus Hulek}
\address{K. Hulek, Institut f\"ur Algebraische Geometrie, Welfengarten 1, 30167 Hannover, Germany}
\email{hulek@math.uni-hannover.de}
\address{Current address: Institute for Advanced Study, School of Mathematics, 1 Einstein Drive, Princeton, NJ 08540, USA}
\email{hulek@ias.edu}

\author[A. Sch\"urmann]{Achill Sch\"urmann}
\address{A. Sch\"urmann, Universit\"at Rostock, Institute of Mathematics, 18051 Rostock, Germany}
\email{achill.schuermann@uni-rostock.de}

\thanks{We thank the referee for useful remarks which helped to improve the presentation of the text.
The first author has been supported by the Croatian Ministry of Science, Education and Sport under contract 098-0982705-2707 and by the Humboldt Foundation.
The second author gratefully acknowledges support by DFG grant Hu 337/6-2 and the Fund for Mathematics to the  Institute of Advanced 
Study in  Princeton, which provided excellent working conditions. }

\title[Perfect form compactifications]{Smoothness and singularities of the perfect form and the second Voronoi compactification of ${\mathcal A}_g$}

\date{}

\maketitle

\begin{abstract}
We study the cones in the first Voronoi or perfect cone decomposition of quadratic forms with  respect to the question which of these cones are basic or 
simplicial.  
As a consequence we deduce that the singular locus of the moduli stack $\Perf$, the toroidal compactification of the moduli space of principally polarized abelian varieties of 
dimension $g$
given by this decomposition, has codimension $10$ if $g \geq 4$. Moreover we describe 
the non-simplicial locus in codimension $10$. We also show that the second Voronoi compactification $\Vor$ has singularities in codimension $3$ for $g\geq 5$.
\end{abstract}

\section{Introduction}

Reduction theory of quadratic forms is a classical topic. This leads to interesting admissible $\GL_g(\ZZ)$-in\-var\-iant tessellations of the rational closure 
of the space of positive-definite real quadratic forms in  dimension $g$ -- also known as admissible rational polyhedral 
decompositions or fans. 
In algebraic geometry these  
give rise to  toroidal compactifications of the moduli space ${\mathcal A}_g$ of 
principally polarized abelian varieties. This theory was first developed by
Ash, Mumford, Rapoport and Tai \cite{AMRT} and has since been used by numerous authors in geometric studies of moduli of abelian varieties, where it is often essential 
to consider not only the non-compact spaces ${\mathcal A}_g$, but to work with good projective models.
More precisely, such a decomposition defines both a variety and a stack (see \cite{FC}, \cite{SB}).
Toroidal compactifications have the property that the boundary is ``big'', i.e. has codimension $1$. All toroidal compactifications map surjectively onto the Satake
compactification~$\Sat$ which is set-theoretically given by
\begin{equation}\label{equ:beta}
  \Sat = {\mathcal A}_g \sqcup {\mathcal A}_{g-1} \sqcup \ldots \sqcup {\mathcal A}_0.
\end{equation}

In the literature three different types of decompositions have been studied in detail:  the {\em first Voronoi} or {\em perfect cone} decomposition,
the  {\em second Voronoi} decomposition and the {\em central cone} decomposition, leading to the corresponding 
toroidal compactifications $\Perf$, $\Vor$ and $\Igu$ correspondingly, see  \cite{Nam} for more details.  In recent years the meaning of these 
various toroidal compactifications has been clarified. The central cone compactification has been identified with
the {\em Igusa} compactification, a blow-up of the Satake compactification $\Sat$.  The second Voronoi compactification $\Vor$ has a meaningful
interpretation in terms of moduli of abelian varieties. Its boundary points correspond to degenerate abelian varieties, more precisely to semi-abelic
varieties. 
For details the reader is referred to work by Alexeev \cite{alexeev} and Olsson \cite{olsson}. Finally, the first Voronoi or perfect cone 
compactification is a good model to work with from the point of view of the classification theory of higher-dimensional algebraic varieties. Shepherd-Barron \cite{SB} has proved that $\Perf$ has canonical 
singularities and that its canonical bundle is ample if $g \geq 12$. Thus $\Perf$  is a {\em canonical model} in the sense of the minimal model
program if $g\geq 12$.

The construction of toroidal compactifications of the moduli space ${\mathcal A}_g$ very roughly works as follows. 
One first has to choose an admissible rational polyhedral decomposition of the rational closure 
of the space of positive definite symmetric $g \times g$
matrices. 
This defines a compactification  $\tor$ of ${\mathcal A}_g$ as a compact analytic space by adding a stratum to each cone in the chosen decomposition 
(where strata of $\GL_g(\ZZ)$-equivalent cones will be identified).
The codimension of the stratum which is added equals the dimension of the cone.  In the cases we have mentioned above the 
compactification $\tor$ is in fact a projective variety. Naturally the geometric properties of $\tor$ depend essentially on the properties of the chosen
decomposition.

In this note we are especially interested in the singularities of toroidal compactifications. 
The singularities of the {\em variety} given by a fan arise in two different ways. First of all the symplectic group
$\operatorname{Sp}(2g,\mathbb Z)$ has torsion (different from $\pm \operatorname {id}$). This gives rise to finite quotient singularities in ${\mathcal A}_g$ 
(unless the torsion element is a reflection, which only happens in genus $g=2$). Similarly such quotient singularities can arise in the boundary 
of a toroidal compactification due to non-neatness of $\operatorname{Sp}(2g,\mathbb Z)$. Such quotient singularities are well behaved from an
algebraic-geometric point of view and can for many 
considerations be neglected: 
if one replaces the group $\operatorname{Sp}(2g,\mathbb Z)$ by a principal congruence subgroup 
$\Gamma(\ell)=\{g \in  \operatorname{Sp}(2g,\mathbb Z) \mid g \equiv {\bf 1} \mod \ell \}$, 
then for $\ell \geq 3$ this group is neat and in particular torsion free. Hence the corresponding level covers  ${\mathcal A}_g(\ell)$ respectively $\torl$ will not
acquire such singularities. A second type of singularities arises from ``bad'' behavior of cones $\sigma$ in the chosen decomposition. Recall that a cone
$\sigma$  whose general element has rank $g$ is called {\em basic} if the reduced integral generators of its $1$-dimensional faces can be completed 
to a $\mathbb Z$-basis of $\operatorname{Sym}^2({\mathbb Z}^g)$. By an {\em integrally reduced generator} we mean that the gcd of its entries equals $1$. Whenever we speak about integral generators we will from now on assume them to be integrally reduced.   
A cone is called {\em simplicial} if its rational generators can be completed to a  $\mathbb Q$-basis, in other words if the generators are linearly independent.
If a cone $\sigma$ is basic, then its corresponding stratum lies in the non-singular locus of $\tor(\ell)$ for $\ell \geq 3$,
whereas simplicial cones give rise to finite quotient singularities by abelian groups, the group being the quotient of the lattice $\operatorname{Sym}^2({\mathbb Z}^g)$ by
the lattice spanned by the integral generators of the cone. Non-simplicial cones give rise to more general singularities.
Taking a level cover will not remove singularities which arise from non-basic cones, these are then singularities of the corresponding {\em stack}. 
We thus obtain  open sub-stacks  $\torb \subset \tors \subset \tor$ given by the partial compactification defined by the basic and simplicial singularities respectively.

Knowledge of the singularities of a toroidal compactification $\tor$ and its level covers is obviously of geometric interest. One reason is that singularities are crucial 
for the  understanding of the birational geometry of these varieties, the other is that they are also important in order to understand topological properties.
One example for the latter is the 
work by Grushevsky, Tommasi and the second author \cite{GHT} on stable cohomology of $\Perf$. The relevance of 
${\mathcal A}_{g,\operatorname{smooth}}^{\operatorname {Perf}}$ (and also ${\mathcal A}_{g,\operatorname{simp}}^{\operatorname {Perf}}$) in this 
context is that Poincar\'e duality 
holds here. 
 For $g\leq 3$ the first and second Voronoi decomposition as well as
the central cone decomposition all coincide and all cones are basic. Hence the stacks are smooth in genus $\leq 3$. This changes in genus $4$.
All cones in  the second Voronoi decomposition are still basic, but this is no longer the case in the first Voronoi or perfect cone decomposition, which in
genus $4$ coincides with the central cone decomposition.  Here, there is one non-basic cone, namely the perfect cone of the root lattice $\mathsf{D}_4$, which has dimension $10$ and $12$ rays, hence is neither basic nor simplicial.
This defines the unique singular point in the $10$-dimensional stack ${\mathcal A}_4^{\operatorname{Perf}}$.
This also means that for $g\geq 4$ the stack $\Perf$ will always be singular in codimension $10$ (or less). 
There is no a priori reason that the codimension of the singular locus of $\Perf$ could not be less than $10$ for $g > 4$. Our aim is to show that this is 
not the case. We shall also classify the non-simplicial locus in codimension~$10$ -- as it turns out, it all comes from the root lattice $\mathsf{D}_4$.

From the point of view of tessellations of the rational
closure of the cone of positive-definite real quadratic forms our results show the following: with respect to the the properties basic or simplicial the perfect cone decomposition
is much more uniform up to dimension $9$
than one might hope for in view of the explosion of arithmetically inequivalent perfect forms from dimension $6$ onwards. 
The only exception to this in dimension $10$ is given by the cone of the root lattice $\mathsf{D}_4$, which emphasizes once more the extraordinary 
importance of this root lattice.
 
\begin{theorem}\label{theo:main}
(i) Every cone of dimension at most $9$ in the perfect cone decomposition is basic, so the  integral generators can be completed to a $\ZZ$-basis of
$\Sym^2(\ZZ^g)$. 

(ii)  With the exception of the cone of the root lattice $\mathsf{D}_4$, every cone in the perfect cone decomposition of dimension at most $10$ is simplicial, so its integral 
generators can be completed to a $\QQ$-basis of $\Sym^2(\QQ^g)$. 
\end{theorem}

Form this we immediately obtain:
\begin{corollary}\label{cor:main}
The stack $\Perf$ is smooth for $g \leq 3$ and
the codimension of both the singular and the non-simplicial substack  of
$\Perf$ is $10$ if $g \geq 4$.
\end{corollary}

Before we give the proof of this theorem  we note the following consequence for 
the intersection cohomology of the variety~$\Perf$:

\begin{corollary}[\cite{GHT}]
In degree $k \leq 10$ the intersection cohomology of the variety $\Perf$ is isomorphic to its singular cohomology:
$$
\operatorname{IH}^k(\Perf) \cong \operatorname{H}^k(\Perf) \quad {\mbox {for}}  \quad k \leq 10.
$$
\end{corollary}

We also answer the same question for the second Voronoi compactification.

\begin{theorem}\label{theo:main2}
(i) For $g\leq4 $ every cone in the second Voronoi compactification is basic.

(ii) For $g\geq 5$ there are  non-simplicial cones in dimension $3$,  thus for these dimensions  $\Vor$ is singular in codimension $3$. 
\end{theorem}

Our results lead to  natural further questions which we will pose as
problems at the end of Sections~\ref{sec:ninedim}, \ref{LowDimConeSimplicial}
and \ref{sec:proof2}.

\section{Proof of Theorem~\ref{theo:main}}   \label{sec:ProofOfThm}

\subsection*{Background on perfect cones}

Let us start with some background on the perfect cone decomposition.
For a vector $v\in \ZZ^g$ we write $p(v)=v v^t$ for the corresponding rank-$1$ matrix.
Any cone $\sigma$ of the perfect cone decomposition 
is of the form $\sum_{i=1}^M \RR_+ p(v_i)$, where $v_1, \ldots , v_M\in \ZZ^g$. 
More precisely, perfect cones arise in this way,
if and only if the $v_i$ are coordinate vectors (with respect to a lattice basis)
of {\em all} the minimal vectors for some $g$-dimensional lattice.
Here, minimal vectors come in pairs $\pm v_i$ and because
of $p(v_i)= p(-v_i)$ it suffices to consider one coordinate vector for each of them.
Moreover, by the minimal vector property the $\gcd$ of the coordinates is~$1$, for each~$v_i$.
So the corresponding generators $p(v_1), \ldots, p(v_M)$ of the perfect cone $\sigma$
are integrally reduced (that is, their $\gcd$ is~$1$ too).
For further details on perfect cones (which are sometimes called Voronoi domains) 
and the associated theory of perfect quadratic forms 
we refer to~\cite{MartinetBook} and \cite{schuermann-2009}.

By $\dim\sigma$ we denote the {\em dimension} of a cone $\sigma$ in
$\operatorname{Sym}^2({\mathbb Z}^g)$, that is, the dimension of its linear
hull.
The faces of perfect cones (intersections with a supporting hyperplane)  
are perfect cones of lower dimensions themselves.
For a given~$g$,
the perfect cones give a face-to-face tessellation of the rational closure
of the space of real, positive-definite $g\times g$ matrices.
This tessellation is invariant with respect to
the action $M\mapsto U M U^t$
of the group $\operatorname{GL}_g(\ZZ)$ on this space.
Two perfect cones are called $\operatorname{GL}_g(\ZZ)$-equivalent,
if they are in the same orbit with respect to this group action.
By a classical theorem due to Voronoi,
we know that there exist only finitely
many perfect cones up to equivalence
for every fixed~$g$.
For a given~$g$ it is therefore possible, at least in principle, to classify all
perfect cones~$\sigma$ of a given dimension~$\dim \sigma = N$.
However, for a classification of all perfect cones of dimension~$N$,
hence in spaces $\operatorname{Sym}^2({\mathbb Z}^g)$ with varying $g$,       
we need an additional argument.

\subsection*{Reduction to finitely many perfect cones}

Let us take a cone $\sigma$ of dimension $N$ of the perfect cone decomposition, generated by $p(v_1), \dots, p(v_M)$.
Denote by $d$~the dimension of the lattice $L=\ZZ v_1 + \dots + \ZZ v_M$ in $\ZZ^g$.
Then $L$ is a finite index sublattice of 
the lattice $L' = (L\otimes \RR)\cap \ZZ^g$ having the same affine span.
$L'$ has a $\ZZ$-basis $\{w_1, \dots, w_d\}$ that can be 
extended to a $\ZZ$-basis $\{w_1, \dots, w_g\}$ of $\ZZ^g$.
Thus by applying a suitable $\operatorname{GL}_g(\ZZ)$ mapping 
(base change for $\ZZ^g$) we may assume that $L'$ is equal to~$\ZZ^d$. 
As a consequence, the properties that concern us, namely that the family $\left(p(v_i)\right)_{1\leq i\leq M}$ is extensible to a $\ZZ$- or $\QQ$-basis, 
can be considered by assuming $d=g$.
I.e., we may assume that $L$ is a full dimensional, finite index sublattice of~$\ZZ^g$.
Since $p(v_1), \ldots, p(v_g)$ are linearly independent if $v_1, \ldots, v_g$ 
are linearly independent, we may assume $g \leq \dim\sigma$. 
This leaves only a finite number of cases to consider,
when classifying all perfect cones of a given dimension~$N$.

\subsection*{Proof of (i)}

Unfortunately we are not aware of a general method for proving basicness of perfect cones without a full enumeration. So for the proof of (i) we have to get a hand on all~$9$ dimensional perfect cones for all $g\leq 9$ and check for each one of the cones 
that they are basic.
Note that given the reduced integral generators $p(v_1), \ldots , p(v_M)$ 
of a fixed perfect cone, 
we can easily check basicness computationally,
for instance, by considering the generators themselves as
vectors and checking if their Gram matrix has determinant~$1$.

We split the needed classification of $9$-dimensional perfect cones 
into three parts:
\begin{itemize}
\item For $g\leq 7$ all perfect cones have been classified up to
$\operatorname{GL}_g(\mathbb{Z})$ equivalence in~\cite{PerfFormModGrp}. 

\item For $g=8$ we computationally classify all
          $9$-dimensional cones in our Lemma~\ref{TheEnum} below. 
 
\item For $g=9$, the perfect cones of dimension~$9$ 
          have to be simplicial as we show in our Lemma~\ref{SimplicialityLemma} 
          in the next section.
          Therefore, these perfect cones are in $1$-to-$1$-correspondence with
          sets of linearly independent 
          vector pairs $\pm v_1,\ldots, \pm v_9 \in \ZZ^9$,
          that are coordinates  (with respect to a lattice basis)
          of {\em all} minimal vectors for some $9$-dimensional lattice.
          Among others, such vector configurations 
          were classified in~\cite{MinkowskianSublattices}, 
          giving us $31$~perfect cones of dimension~$9$ for $g=9$,
          corresponding to the entries with
          $s=r=s'=9$ of Tables~2,~3 and~7 in~\cite{MinkowskianSublattices}.

\end{itemize}

In addition to the summarized information given here, 
we provide a complete list of all $9$-dimensional perfect cones in electronic form 
on the webpage~\cite{MathieusWebpage}. Representatives of the cones $\sigma$ of
the perfect decomposition of dimension at most $10$ and genus $g\leq 9$ are given
there. The classification in genus $g=8, 9$ and dimension $10$ was obtained by
extensions of the method of Section \ref{sec:ninedim}.
Table \ref{NrOrbitCones} gives the number of orbits of cones.

\begin{table}
\caption{Number of orbits of cones in the perfect cone decomposition for $g\leq 9$ and dimension at most $10$}
\label{NrOrbitCones}
\begin{tabular}{||c||c|c|c|c|c|c|c||}
\hline
\hline
g$\downarrow$, dim$\rightarrow$ 
          & 4 & 5 & 6 & 7  & 8  & 9   & 10\\
\hline
\hline
4         & 1 & 3 & 4 & 4  & 2  & 2   & 2\\
\hline
5         &   & 2 & 5 & 10 & 16 & 23  & 25\\
\hline
6         &   &   & 3 & 10 & 28 & 71  & 162\\
\hline
7         &   &   &   & 6  & 28 & 115 & 467\\
\hline
8         &   &   &   &    & 13 & 106 & 783\\
\hline
9         &   &   &   &    &    & 44  & 759\\
\hline
\hline
\end{tabular}
\end{table}

\subsection*{Proof of (ii)}

For our proof of assertion (ii) we ``only'' have to consider the~$10$-dimensional perfect cones, because from (i) we know that all cones of dimension $\leq 9$ 
are simplicial (as they are basic, which is stronger).
It seems quite challenging to classify all $10$-dimensional perfect cones, so our
proof relies on a combination of mathematical reasoning and computer assisted checks:
\begin{itemize}
\item 
For $g\leq 7$ we can use the previously known classification of perfect cones again.
According to \cite{PerfFormModGrp} there are $656$ inequivalent
$10$-dimensional perfect cones
and we checked computationally that all of them are simplicial, except the one of the root lattice~$\mathsf{D}_4$.

On the webpage \cite{MathieusWebpage} we provide a complete list of all 
$10$-dimensional perfect cones with $g\leq 7$.

\item
For $g=8$, $g=9$ and $g=10$ 
there is no complete classification of $10$-dimensional perfect cones known so far.
Here, we use our Lemma~\ref{SimplicialityLemma} below. It shows 
for general $g$, that all perfect cones of dimension $g$, $g+1$ and $g+2$ are necessarily simplicial.

\end{itemize}

\qed

\subsection*{Classifying $10$-dimensional perfect cones?}

A complete classification of $10$-dimensional perfect cones
appears to be highly challenging, if not out of reach at the moment.
Although an extension of our method to prove Lemma~\ref{TheEnum} could possibly 
be feasible (to deal with the case $g=9$), the necessary extension of the
work in \cite{MinkowskianSublattices} to the case $g=10$ seems 
computationally hardly realizable.
A slightly simpler task could be the case $g=8$. 
More generally, enumerating the number of inequivalent 
$(8+k)$-dimensional perfect cones for small~$k$ 
could be feasible for~$g=8$. 
Lemma~\ref{TheEnum} shows that there are $106$ classes for $k=1$.
Since this number is much less than the $10916$ classes of full dimensional perfect cones for $g=8$, which have been classified in~\cite{perfectdim8},
one may hope to obtain a complete list for $k=2$ as well. 

\medskip

Note that only with the help of our Lemma~\ref{SimplicialityLemma} below
we could avoid this very difficult (currently impossible) classification of~$10$-dimensional perfect cones, in order to prove part (ii) of the Theorem.
On the other hand, for part (i) of the Theorem, our proof 
relies on the complete classification of all $9$-dimensional perfect cones, 
which would have been extremely painful (if not impossible) 
without computer assistance.
So, it is just this special combination of mathematical reasoning and
computer assistance which 
allowed us to obtain the results of this paper.

\section{Low dimensional perfect cones are simplicial}\label{LowDimConeSimplicial}

Perfect cones for $g=2$ and $g=3$ are known to be simplicial.
In fact, in these two cases there exists only one class of
top-dimensional perfect cones of dimension $3$ and $6$ respectively.
Both are associated to the root lattice $\mathsf{A}_g$
and are known to be simplicial.
Since all perfect cones are faces of top-dimensional perfect cones
and since faces of simplicial cones are simplicial, we see that 
perfect cones of any dimension~$N$ are simplicial for $g=2$ and $g=3$.
This argument does not work for $g\geq 4$ 
as, for instance, the top-dimensional perfect cones associated to the root lattices $\mathsf{D}_g$
are not simplicial in these cases.
The following Lemma shows, however, that perfect cones are simplicial 
for any $g$ in case their dimension is sufficiently small, 
namely at most~$g+2$.

\begin{lemma}\label{SimplicialityLemma}
Let $\sigma=\RR_+ p(v_1) + \dots \RR_+ p(v_M)$ be a cone of
the perfect cone decomposition of $\operatorname{Sym}^2({\mathbb Z}^g)$.
Assume that $(v_1,\dots, v_M)$ span $\RR^g$ and that $\dim\sigma=g$,
$g+1$ or $g+2$.
Then $\sigma$ is simplicial.
\end{lemma}

\proof 
As mentioned above, the cases $g=2$ and $g=3$ are known by classical 
results. So we may assume that $g\geq 4$ and that the 
Lemma is true for all dimensions less than $g$.

Among the generators of $\sigma$ we can find 
$g$ linearly independent vectors $(v_i)_{1\leq i\leq g}$
that determine a basis of $\QQ^g$.
By using this basis, the vector space spanned by $(p(v_i))_{1\leq i\leq g}$
can be identified with the space of diagonal $g\times g$ matrices.
That is, after a suitable base change, 
we may assume that the $v_i$ are the standard basis vectors.
Given an arbitrary vector $v=\sum_{i=1}^g \alpha_i v_i$, we see that  
$p(v)$ is linearly independent of the generators $p(v_i)$, if and only if
$\alpha_i\not=0$ for at least two indices $i_1 \not= i_2$,
since the non-diagonal entry $(i_1, i_2)$ is non-zero in that case.

\medskip

Suppose $\dim \sigma = g$. 
Then any additional generator $p(v)$ of a perfect cone
would belong to the vector space spanned by the $(p(v_i))_{1\leq i\leq d}$.
Therefore, by the argument above, 
$v=\sum_{i=1}^g \alpha_i v_i$ with $\alpha_i\not= 0$ for just one $\alpha_i$.
Hence, $v$ is a multiple of one $v_i$ in that case,
and therefore $p(v)$ is a multiple of $p(v_i)$ as well,
showing that $p(v)$ is no additional generator.

\medskip

Suppose $\dim \sigma=g+1$. Then we can find 
a vector $v$ in $L= \ZZ v_{1} + \ldots  \ZZ v_{g} $ 
with $p(v)$ a generator of $\sigma$ being linearly independent of the $p(v_i)$.
After suitable scaling of the vectors $v_i$ we may assume without loss of generality 
that $v=\sum_{i=1}^r v_i$ with $2\leq r\leq g$.
Suppose now that $w=\sum_{i=1}^g \alpha_i v_i$ with $p(w)\in \sigma$. 
Then there exist $\gamma$ and $\beta_i$ such that 
$p(w)=\gamma p(v) + \sum_{i=1}^g \beta_i p(v_i)$
and we obtain the equation
\begin{equation} \label{eqn:quadratic-relation}
\left( \sum_{i=1}^g \alpha_i x_i\right)^2 
= \gamma \left(\sum_{i=1}^r x_i\right)^2  + \sum_{j=1}^g \beta_i x_i^2
\end{equation}
for the corresponding quadratic forms.
From this it follows that we cannot have 
$\alpha_i\alpha_j\not=0$ for some $i\geq r+1$ and $i\not= j$. 
This implies that either 
(i) $\alpha_{i_0}\not=0$ only for one $i_0\geq r+1$ or 
(ii) $\alpha_i=0$ for all $i\geq r+1$.

\smallskip

(i) In the first case we deduce  $w=\alpha_{i_0} v_{i_0}$ from~\eqref{eqn:quadratic-relation}.
Hence $p(w)$ is a multiple of $p(v_{i_0})$ in that case
and therefore $p(w)$ cannot be an additional generator of~$\sigma$.

\smallskip

(ii) In the second case we deduce from~\eqref{eqn:quadratic-relation} that
$\alpha_i\alpha_j = \gamma$ for all $i,j\leq r$ with $i\not=j$.
In particular, $\alpha_i\not= 0$ for all $i\leq r$.

In the case $r=2$ we have $\beta_i=0$ for $i\geq 3$ because of 
$\alpha_i=0$ for all $i\geq 3$. 
Thus the putative generator $p(w)$ is actually in the 
three dimensional vector space spanned by the generators $p(v)$, $p(v_1)$, $p(v_2)$.
We know from the already solved case of dimension~$g=2$
that $p(w)$ must be a multiple of one of the other three
generators, if they span a perfect cone themselves.
Every face of a perfect cone is a perfect cone, so it remains to verify that
$p(v)$, $p(v_1)$, $p(v_2)$ span a face of the perfect cone~$\sigma$.
That is, we have to check whether or not there exists a linear form on 
$\operatorname{Sym}^2({\mathbb Z}^g)$ that evaluates to~$0$ 
on the three generators and is positive on the remaining ones.
Taking the inner product with $p(v_3) + \ldots + p(v_g)$ gives such a linear form.

So suppose $r\geq 3$. 
Then at least two of the $\alpha_i$ with $i\leq r$ 
are of the same sign and so $\gamma >0$.
So, all $\alpha_i$ with $i\leq r$ are of the same sign, which 
we can assume to be positive.
For all subsets $S=\{i_1, i_2, i_3\}$ of $\{1,\dots,r\}$ 
with distinct indices $i_1, i_2, i_3$ 
the equations $\alpha_{i_1}\alpha_{i_2}=\alpha_{i_1}\alpha_{i_3}=\alpha_{i_2}\alpha_{i_3}=\gamma$ have the unique solution $\alpha_{i_1}=\alpha_{i_2}=\alpha_{i_3}=\sqrt{\gamma}$.
So $w$ is actually a multiple of $v$, showing again that 
$p(w)$ cannot be an additional generator of $\sigma$.

\medskip

Suppose $\dim \sigma=g+2$. By a suitable scaling of the $v_i$ we can find
vectors of the form $v=\sum_{i=1}^r v_i$ and $v'=\sum_{i=1}^g \lambda_i v_i$ 
such that the set $\{p(v_1), \dots, p(v_g),p(v), p(v')\}$ consists of generators of $\sigma$ 
which form a basis of the rational vector space spanned by $\sigma$.
Assume $p(w)$ is an element of this vector space and $w=\sum_{i=1}^g \mu_i v_i$.
Then there exist $\alpha_i$, $\beta$, $\gamma$ such that in terms
of quadratic forms we have
\begin{equation}  \label{eqn:quadratic-relation2}
\sum_{i=1}^g \alpha_i x_i^2 + \beta \left(\sum_{i=1}^r x_i\right)^2 
+ \gamma \left(\sum_{i=1}^g \lambda_i x_i\right)^2 = 
\left(\sum_{i=1}^g \mu_i x_i\right)^2.
\end{equation}
By the arguments for the case $\dim \sigma=g+1$
we can assume $\beta\not=0$ and $\gamma\not=0$.
We distinguish the two cases: (i) $r<g$ and (ii) $r=g$.

\smallskip

(i) Assume $r<g$. 
Then \eqref{eqn:quadratic-relation2} implies
$\gamma \lambda_i\lambda_g = \mu_i\mu_g$ for all $i<g$. 

If $\lambda_g=0$ then $\mu_i \mu_g=0$.
If $\mu_g\not=0$ then $\mu_i=0$ for all $i<g$ and $w$ is a multiple of~$v_g$.
If $\mu_g=0$ then we have necessarily $\alpha_g=0$ and the problem is reduced to a lower dimensional one.

If $\lambda_g\not=0$, then \eqref{eqn:quadratic-relation2} implies $\mu_g\not=0$.
So $\lambda_i = \mu_i \kappa$ with $\kappa=\mu_g/(\gamma \lambda_g)$ for
$i<g$. 
Hence the quadratic forms $p(v')$ and $p(w)$ are multiples of each other,
showing that $p(w)$ cannot be an additional generator of $\sigma$.

\smallskip

(ii) Assume $r=g$.
Then from \eqref{eqn:quadratic-relation2} we deduce 
the equality $\beta + \gamma \lambda_i\lambda_j = \mu_i \mu_j$
for all $i\not= j$. 
We may assume that all $\lambda_i$ and $\mu_i$ are non-zero, 
since otherwise we can permutate between $v$, $v'$ and $w$ 
and reduce to the preceding case (i). 
By a suitable scaling of the $v_i$ we may assume $\mu_g=1$.
So we get $\mu_i = \beta + \gamma \lambda_i \lambda_g$ and for $1\leq i<j < g$:
\begin{equation*}
\begin{array}{rcl}
0 &=& \mu_i\mu_j - \beta - \gamma \lambda_i \lambda_j\\
&=& (\beta + \gamma \lambda_i\lambda_g)(\beta + \gamma\lambda_j\lambda_g) - \beta - \gamma \lambda_i\lambda_j\\
&=& \beta^2 - \beta + \beta \gamma \lambda_g (\lambda_i + \lambda_j) + \gamma(\gamma \lambda_g^2 -1) \lambda_i \lambda_j\\
&=& a + b(\lambda_i + \lambda_j) + c\lambda_i \lambda_j.
\end{array}
\end{equation*}

If $ac\not= b^2$ then there exists a fractional function $\phi:\PP^1(\RR)\rightarrow \PP^1(\RR)$ such that $\lambda_i = \phi(\lambda_j)$ and $\phi\circ \phi=\operatorname{Id}$.
It follows that there are one or two possible values for $\lambda_i$,
and in case of two different values, we find for any 
pair $(i,j)$ with $1\leq i < j < g$ that $\lambda_i$ is equal to
one value and $\lambda_j$ is equal to the other.
The restriction $g\geq 4$ implies that $\lambda_i$ takes only one value in this case. 

If $ac=b^2$ then we have $c\not=0$ and $(\lambda_i + b/c)(\lambda_j + b/c)=0$.
So, no two $\lambda_i$ with $i<g$ can be different from $-b/c$.
If there is one index $i_0<g $ with $\lambda_{i_0}\not = -b/c$
and if $\lambda_g \not = -b/c$ as well, we may 
choose a different special index than $g$ and obtain a contradiction too.
Otherwise, if $\lambda_g = -b/c$,
we may choose $i_0$ as this special index.
So we may assume $\lambda_i= -b/c$ for all $i<g$.
In case all the $\lambda_i$ are equal (hence also $\lambda_g = -b/c$),
$v$ and $v'$ are multiples of each other, contradicting the
assumption that $p(v)$ and $p(v')$ are linearly independent.

So, we have $\lambda_g \not= -b/c$, $\lambda_i=\lambda_1$ and $\mu_i=\mu_1$ for $i<g$.
As we assume $g\geq 4$, 
setting $\sum_{i=1}^{g-1}x_{i}=0$ in Equation~\eqref{eqn:quadratic-relation2} 
yields $\alpha_i=0$ for $i<g$. 
Using the two variables $y_1=\sum_{i=1}^{g-1} x_i$ and $y_2=x_g$,
we see that Equation~\eqref{eqn:quadratic-relation2} is actually a relation between forms of two variables.
Therefore the putative generator $p(w)$ of $\sigma$ is in 
the three dimensional vector space spanned by the generators $p(v_g)$, $p(v)$, $p(v')$.
We know from the already solved case of dimension~$g=2$
that $p(w)$ must be a multiple of one of the other three
generators, if these span a perfect cone themselves.
As in the proof for the case $\dim \sigma=g+1$ (subcase (ii) with $r=2$),
this follows if we show that $p(v_g)$, $p(v)$, $p(v')$ span a face of the perfect cone~$\sigma$.
To see this we may construct a linear form on $\operatorname{Sym}^2({\mathbb Z}^g)$
that evaluates to~$0$ on the three generators and is positive on the remaining ones.
Such a form is for instance obtained by considering the
inner product with the projection of $p(v_1) + \ldots + p(v_{g-1})$
onto the linear space orthogonal to the span of the three generators $p(v_g)$, $p(v)$ and $p(v')$. \qed

Based on the enumeration results that we have obtained, it seems plausible that simpliciality
also occurs in cases $\dim\sigma=g+3$ or $g+4$.

\section{Classification of $9$-dimensional perfect cones}\label{sec:ninedim}

As mentioned in Section~\ref{sec:ProofOfThm}, for a complete classification of 
$9$-dimensional perfect cones, we need to obtain the full list of all such cones 
for the case~$g=8$.

\begin{lemma}\label{TheEnum}
There are $106$ orbits of $9$-dimensional cones 
in the perfect cone decomposition for $g=8$.
\end{lemma}

\proof
Our classification is based on a computer assisted case distinction.
The full list can be obtained from the webpage~\cite{MathieusWebpage}.
Here, we briefly describe the necessary ingredients and our computational steps:

\subsection*{Using our simplicity lemma}

By Lemma \ref{SimplicialityLemma} 
all $9$-dimensional perfect cones for $g=8$ 
are necessarily simplicial,
 i.e. of the form $\RR_+ p(v_1) + \dots + \RR_+ p(v_9)$, with linearly independent 
generators $p(v_1), \ldots,  p(v_9)$.
As we assume $g=8$, without loss of generality 
we may assume that $(v_1, \dots, v_8)$ are linearly independent too.

Knowing a priori that the $9$-dimensional perfect cones are simplicial
is a quite strong condition. It in particular implies that all its facets
(codimension-$1$-faces) are simplicial $8$-dimensional perfect cones.
Thus we can assume that $\pm v_1,\ldots , \pm v_8 \in \ZZ^8$ are
coordinates (with respect to a lattice basis) 
of all minimal vectors for some $8$-dimensional lattice.
Such sets of eight minimal vector pairs have been classified 
by Martinet in \cite{Martinet} (see entries satisfying $n=s=r=s'=8$ in Tableau 11.1).
Up to $\GL_8(\ZZ)$ equivalence there are $13$~such sets of vectors
and we can assume that $(v_1, \dots, v_8)$ is a representative for one of 
these $13$~orbits.

\subsection*{Reduction to a finite enumeration}

If $L$ denotes the sublattice spanned by the vectors $(v_1, \dots, v_8)$ in $\ZZ^8$,
then from Martinet's classification~\cite{Martinet} we know that the
index $i(L)$ in $\ZZ^8$ satisfies $1\leq i(L)\leq 5$.
We may assume that $(v_1, \dots, v_8)$ is of maximal index 
among all possible $8$-subsets of $(v_1, \dots, v_9)$.

Given a fixed set of vectors $(v_1, \dots, v_8)$, we derive 
some restrictive conditions for the additional vectors $v_9\in\ZZ^8$ 
from this maximality assumption: 
For $j=1,\ldots, 8$ we define the lattice
\begin{equation*}
L_j=\ZZ v_1 + \dots + \ZZ v_{j-1} + \ZZ v_{j+1} + \dots + \ZZ v_8 + \ZZ v_9.
\end{equation*}
If $L_j$ is full-dimensional then its index $i(L_j)$ in $\ZZ^8$ 
is at most $i(L)$ by our maximality assumption for $L$.
Thus
\begin{equation}  \label{eqn:DetCond}
\left\vert \det(v_1, \dots, v_{j-1}, v_{j+1}, \dots v_8, v_9)\right\vert \in \{0,1, \dots, i(L)\}
,
\end{equation}
with the determinant being~$0$ if and only if $L_j$ is not full-dimensional.
By using the exterior product, we can rewrite these determinant  
conditions as $-i(L) \leq \langle w_j, v_9\rangle \leq i(L)$ 
with $w_j$ being a fixed integral vector,
depending only on the vectors $v_1, \dots, v_{j-1}, v_{j+1}, \dots v_8$
(orthogonal to it).
Thus for a fixed set of vectors $(v_1, \dots, v_8)$ we obtain 
eight pairs of linear inequalities for the possible additional integral vector~$v_9$.
Geometrically, these conditions 
define a parallelepiped and we 
can use for instance the program {\tt zsolve} 
from \cite{4ti2} to enumerate its integral points.

Note that the eight determinant conditions \eqref{eqn:DetCond}
have altogether $(1+ 2i(L))^8$ possible values. It is easy to see
that the matrix $(w_j)_{j=1,\ldots,8}$ is the adjugate matrix of 
$(v_i)_{i=1,\ldots,8}$ and so is of determinant $i(L)^7$.
This implies that the number of integral solutions $v_9$ is actually
reasonable, 
allowing {\tt zsolve} to find all solutions.
So the condition that the lattice spanned by $\{v_1, \dots, v_8\}$ 
is of maximal index is a key assumption for the enumeration to work.

\subsection*{Exploiting symmetry}

We are left with a finite number of possible sets $\{v_1, \dots, v_9\}$.
Using Algorithm~1 of~\cite{MinkowskianSublattices},
we can test for each such set wether it defines a $9$-dimensional 
perfect cone or not.
However, each of these tests is computationally very expensive and 
it is therefore advisable to consider further reductions beforehand, 
in order to finish the classification.

A large reduction of cases can be obtained from using the symmetry 
within a given configuration ${\mathcal V} = (\pm v_1, \dots, \pm v_8)$: 
We only need to 
consider vectors~$v_9$ up to the automorphism group (a subgroup of $\GL_8(\ZZ)$) 
of the configuration ${\mathcal V}$. 
Once we have this integral automorphism group of~${\mathcal V}$
we identify orbits of possible extensions $v_9$ and choose only one 
representative for each orbit to proceed.

To compute the automorphism group, we first consider
the group $G_2$ of rational automorphisms of ${\mathcal V}$, 
which is the hyperoctahedral group
iof size $2^8 8!$. 
In order to obtain the integral automorphisms of ${\mathcal V}$, 
we first determine a subgroup $G_1$ formed by transpositions $(i,j)$ 
and sign changes $(v_i\mapsto -v_i)$ which induce integral automorphisms. 
The full automorphism group within $\GL_8(\ZZ)$ 
is then obtained by iterating over double cosets $G_1 h G_1$ 
and keeping the ones that preserve ${\mathcal V}$.

\subsection*{Treatment of the remaining tuples}

For each set $\{v_1, \dots, v_9\}$ that remains we 
check if each of its $8$-subsets is integrally equivalent  
to one of the $13$~orbits classified by Martinet in~\cite{Martinet}.
This is a necessary condition, as every facet of the potential perfect cone
has to be a perfect cone itself.
If all these checks are positive, we 
test if the candidate $\{v_1, \dots, v_9\}$ defines a $9$-dimensional 
perfect cone by using Algorithm~1 of \cite{MinkowskianSublattices}.

We thus get $131$ systems of $9$ vectors to test for unimodular
equivalence. 
This gives $106$ systems, respectively orbits, under $\GL_8(\ZZ)$-equi\-va\-lence.

\medskip

The computation was dominated by the realizability tests using 
Algorithm~1 of \cite{MinkowskianSublattices}.

\qed

\subsection*{Open problem}

There are examples of cones which are simplicial but not basic. 
One such example is given by
the shortest vectors of the dual root lattice $\mathsf{E}_7^*$ (the index of the corresponding sublattice of $\operatorname{Sym}^2({\mathbb Z}^7)$ is $384$).
Our classification of $9$-dimensional perfect cones shows that
perfect cones are simplicial and basic in low dimensions.
Note however that this is not the case for general polyhedral cones. 
If we take $v_1=(1,1)$ and $v_2=(1,-1)$ for instance, then the relation
\begin{equation*}
\frac{1}{2} \left( p(v_1) + p(v_2) \right) = p((1,0)) + p((0,1))
\end{equation*}
shows  that the cone spanned by $\{p(v_1), p(v_2)\}$ is simplicial but not basic,
since the two generators cannot be extended to a $\ZZ$-basis of $\operatorname{Sym}^2({\mathbb Z}^2)$. 
Of course, $\{\pm v_1, \pm v_2\}$ cannot be realized as the set of shortest vector of a $2$-dimensional lattice, so $\{p(v_1), p(v_2)\}$ does not generate a perfect cone.
In contrast, our classification shows for perfect cones and $g\leq 9$:
If $\{v_1, \dots, v_g\}$  are independent vectors of ${\mathbb Z}^g$, such that $\{p(v_1),\dots, p(v_g)\}$ are the generators of a perfect cone,
then the cone is simplicial and basic.
Is this true for general~$g$?

\section{Proof of Theorem~\ref{theo:main2}}\label{sec:proof2}

Here we consider the singularities of the second Voronoi compactification $\Vor$.
The second Voronoi decomposition is described in detail in \cite{schuermann-2009}.
It is known only up to dimension~$5$ (see \cite{engel}).
For $g\leq 5$, a complete system of $\GL_g(\ZZ)$-inequivalent polyhedral cones in~$\operatorname{Sym}^2({\mathbb Z}^g)$
can electronically be accessed using the program~\cite{scc} for instance.

For the proof of Theorem~\ref{theo:main2} (i) one can simply check whether or not 
all polyhedral cones for $g\leq 4$ are basic. For this it is even enough to check
if the top dimensional cones of the second Voronoi decomposition are basic.
There is only one such cone each for $g=2$ and $g=3$ and there are three cones
for $g=4$, as already observed by Voronoi (see~\cite{voronoi}). 
We do not know if Voronoi checked them for basicness, but he could certainly have done so.
Recall that basicness of a cone follows if the Gram matrix of its generators has determinant~$1$.

All cones of dimension $1$ or $2$ are spanned by one or two extremal rays and are thus trivially simplicial.
Thus for the proof of statement (ii) of Theorem~\ref{theo:main2}, it suffices 
to find a non-simplicial cone of dimension~$3$ for $g\geq 5$.
Using \cite{engel} we find two cones of dimension~$3$ that are spanned by $4$ generators and so are non-simplicial. They are:
\begin{equation*}
\begin{array}{c}
\left(\begin{array}{ccccc}
4&-2&-2&0&-2\\
-2&4&0&-1&1\\
-2&0&4&-1&1\\
0&-1&-1&3&-1\\
-2&1&1&-1&3
\end{array}\right)
,\mbox{~}\left(\begin{array}{ccccc}
2&-1&-1&0&-1\\
-1&2&0&0&0\\
-1&0&2&-1&1\\
0&0&-1&2&-1\\
-1&0&1&-1&2
\end{array}\right)
,\mbox{~}\\
\left(\begin{array}{ccccc}
2&-1&-1&0&-1\\
-1&2&0&-1&1\\
-1&0&2&0&0\\
0&-1&0&2&-1\\
-1&1&0&-1&2
\end{array}\right)
,\mbox{~}\left(\begin{array}{ccccc}
0&0&0&0&0\\
0&0&0&0&0\\
0&0&0&0&0\\
0&0&0&1&-1\\
0&0&0&-1&1
\end{array}\right)
\end{array}
\end{equation*}
and
\begin{equation*}
\begin{array}{c}
\left(\begin{array}{ccccc}
3&-1&-1&-1&-1\\
-1&4&-1&-1&0\\
-1&-1&3&1&-1\\
-1&-1&1&3&-1\\
-1&0&-1&-1&4
\end{array}\right)
,\mbox{~}\left(\begin{array}{ccccc}
6&-2&-2&-2&-2\\
-2&6&-1&-1&0\\
-2&-1&4&1&-1\\
-2&-1&1&4&-1\\
-2&0&-1&-1&6
\end{array}\right)
,\mbox{~}\\
\left(\begin{array}{ccccc}
2&0&-1&-1&-1\\
0&2&-1&-1&0\\
-1&-1&2&1&0\\
-1&-1&1&2&0\\
-1&0&0&0&2
\end{array}\right)
,\mbox{~}\left(\begin{array}{ccccc}
5&-1&-2&-2&-2\\
-1&4&-1&-1&0\\
-2&-1&3&1&0\\
-2&-1&1&3&0\\
-2&0&0&0&4
\end{array}\right)
\end{array}
\end{equation*}
The two polyhedral cones spanned by these generators 
can be extended to spaces $\operatorname{Sym}^2({\mathbb Z}^g)$
for any $g>5$ and so the result follows. 

\qed

\subsection*{Open problem}

We note that we did not find any polyhedral cone of index $>1$ in the
second Voronoi decomposition.
It appears to be an open problem whether or not for every~$g$, 
all cones in the second Voronoi decomposition have generators 
of its $1$-dimensional faces that generate an integral lattice in 
$\operatorname{Sym}^2({\mathbb Z}^g)$ of index~$1$.

\end{document}